\documentclass[11pt]{article}
\usepackage{amsfonts}
\usepackage{epsfig}

\usepackage{amssymb}
\usepackage{subfigure}
\usepackage{amsmath}
\usepackage{float}
\usepackage{tabularx}
\usepackage{graphicx}
\usepackage{verbatim}

\usepackage{amsthm}
\usepackage{amsfonts}
\usepackage[T1]{fontenc}
\usepackage[utf8x]{inputenc}

\usepackage[lflt]{floatflt}
\usepackage{graphicx}
\usepackage{amsmath}
\usepackage{amssymb}
\usepackage{fancyhdr}
\usepackage{color}
\usepackage[affil-sl]{authblk}
\usepackage{cancel}
\usepackage{subfigure}
\usepackage{fancyhdr}
\usepackage{makeidx}
\usepackage{latexsym}
\usepackage{hyperref}

\usepackage{color}
\usepackage{subfigure}
\usepackage{verbatim}
\usepackage{enumerate}
\usepackage[usenames,dvipsnames]{xcolor}

\newtheorem{teo}{Theorem}

\newtheorem{defi}{Definition}
\newtheorem{propo}{Proposition}
\newtheorem{lema}{Lemma}

\newtheorem{obs}{Remark}

\makeindex \textwidth=17cm \textheight=25cm \oddsidemargin=-0.5cm
\newcommand{\CC}{{\mathcal C}}
\newcommand{\bbR}{{\mathbb R}}

\newcommand{\bbC}{{\mathbb C}}

\newcommand{\e}{\epsilon}
\newcommand{\al}{\alpha}
\newcommand{\la}{\lambda}
\newcommand{\p}{\partial}
\newcommand{\be}{\beta}
\newcommand{\G}{\Gamma}

\begin{document}
\begin{center}\emph{}
\LARGE
\textbf{An integral Relationship for a new Fractional One-phase Stefan  Problem}
\end{center}

                   \begin{center}
                  {\sc Sabrina D. Roscani and Domingo A. Tarzia}\\
 CONICET - Depto. Matem\'atica,
FCE, Univ. Austral,\\
 Paraguay 1950, S2000FZF Rosario, Argentina \\
(sabrinaroscani@gmail.com, dtarzia@austral.edu.ar)
                   \vspace{0.2cm}

       \end{center}
      
\small

\noindent \textbf{Abstract: }

  In this paper a new one-dimensional fractional one-phase Stefan problem with a temperature boundary condition at the fixed face is considered. An integral relationship between the temperature and the free boundary is obtained which is equivalent to the fractional Stefan condition. Moreover, an exact  solution of similarity type expressed in terms of Wright functions is given.


\noindent \textbf{AMS:} Primary: 35R35, 26A33, 35C05. Secondary: 33E20, 80A22.  \\

\noindent \textbf{Note:} This paper is now published (in revised form) in Fract. Calc. Appl. Anal., Vol. 21, No 4 (2018), pp. 901-918, DOI:10.1515/fca-2018-0049, and is available online at \url{http://www.degruyter.com/view/j/fca}, so always cite it with the journal's coordinates.


\section{Introduction}
\label{sec:intro}

\noindent The free boundary problems for the one-dimensional diffusion equation  are problems linked to the processes of melting and freezing. In these problems the diffusion, considered as a heat flow, is expressed in terms of instantaneous local flow of temperature, and a latent heat-type condition at the interface connecting the velocity of the free boundary and the heat flux of the temperatures in both phases. This kind of problems have been widely studied in the last 50 years (see \cite{Alexiades, Cannon, Crank, Elliott, Gupta, Lunardini, Rubi, Stefan, Tar1, Tar2}).\\

\noindent In this paper a new fractional Stefan problem is considered. That is, a problem governed by a fractional diffusion equation with two unknown functions: a two variables function $u=u(x,t)$ and a free boundary $s=s(t)$.\\

\noindent Voller et al. \cite{Voller:2013} state that if we consider an ideal non local flow and we replace it in a heat balance equation,  we derive in subdiffusion, modeled by a fractional diffusion equation involving the fractional Caputo derivative in time of order $\al \in (0,1)$. This kind of diffusion is called anomalous diffusion and it has been studied by numerous authors \cite{ FM-AnaPropAndAplOfTheW-Func, Gusev, Kilbas, Luchko3,FM-TheFundamentalSolution,  MK:2000, Podlubny, Yuste:2010}. Several  applications were considered and, in particular, Mainardi studied in \cite{FM-libro}  the application to the theory of linear viscoelasticity.\\

\noindent Regarding the formulation of fractional free boundary problems, it has been observed that the movement of the moisture boundary   in an horizontal porous brick diffusion can be both as  sub or superdiffusive. Other experiments have shown that growth of frost on a cooled plate can be superdiffusive \cite{TBR:1993}. 

\noindent We propose here to formulate our problem by changing the classical Fourier law which sets that the heat flux  is proportional to the temperature gradient. We suppose that the heat flux is proportional to a ``fractional back in time gradient''
\begin{equation}\label{grad fr}  
q=-k^{RL}_{0}D^{1-\al}_t \frac{\p}{\p x} u (x,t),
\end{equation}

\noindent in which the  fractional derivative is the Riemann--Liouville derivative respect on time of order $1-\al$ ($\al $ $\in (0,1)$).  Then, the classical diffusion equation becomes in our governing fractional diffusion equation (hereinafter FDE) given by:  

\begin{equation}\label{FDE}   \frac{\p}{\p t}u(x,t)=\lambda^2\, \frac{\p}{\p x}\left(^{RL}_{0}D^{1-\al}_t \frac{\p}{\p x} u (x,t)\right),
\end{equation}

\noindent where the Riemann--Liouville fractional derivative is given by
$$ ^{RL}_{0}D^{1-\al}_t \frac{\p}{\p x} u (x,t)= \frac{1}{\Gamma(\al)} \frac{\p}{\p t}\displaystyle\int^{t}_{0}\frac{\frac{\p}{\p x} u(x,\tau)}{(t-\tau)^{1-\al}}d\tau, \qquad \al \in (0,1)$$

\noindent and $\G$ is the Gamma function defined by  $\G(x)=\int_0^\infty  w^{x-1}e^{-w}dw$. 

\noindent The order of derivation $1-\alpha$ is linked to the strong relationship between this equation and the fractional diffusion equation associated to the fractional derivative in the Caputo sense.
 
\noindent We will ask to the interface velocity of the free boundary to be  proportional to the fractional back in time gradient (\ref{grad fr}). So, the classical Stefan condition will be replaced by the following ``fractional Stefan condition'':
\begin{equation}\label{cond fracc de stefan}
\frac{d}{dt}s(t)=-\left.^{RL}_{0}D^{1-\al}_t \frac{\p}{\p x} u (x,t)\right|_{(s(t),t)}, \quad 0<t\leq T.
\end{equation}

\noindent Therefore, our interest problem is given by:
 
\begin{equation}{\label{St}}
\begin{array}{llll}
     (i)  &   \frac{\p}{\p t}u(x,t)=\lambda^2\, \frac{\p}{\p x}\left(^{RL}_{0}D^{1-\al}_t \frac{\p}{\p x} u (x,t)\right), &   0<x<s(t), \,  0<t<T, \,  0<\al<1 , \, \,\\
     (ii) &   u(x,0)=f(x), & 0\leq x\leq b=s(0),\\ 
       (iii)  &  u(0,t)=g(t),   &  0<t\leq T,  \\

         (iv) & u(s(t),t)=0, & 0<t\leq T, \, C\geq 0,\\

        (v) & \frac{d}{dt}s(t)=-\left.^{RL}_{0}D^{1-\al}_t \frac{\p}{\p x} u (x,t)\right|_{(s(t),t)}, & 0<t\leq T, 
                                             \end{array}
                                             \end{equation}

\noindent where $f$  and  $g$ are non-negative continuous functions defined in $(0,T]$.\\

\noindent We present in Section 2 some basic concepts of fractional calculus that will be used later. In Section 3 an integral relationship between a solution of problem (\ref{St}),  given by functions $u=u(x,t)$ and $s=s(t)$, is be proved by using the Green's Theorem. We also  prove that, under certain hypothesis, this integral relationship is equivalent to the fractional Stefan condition (\ref{cond fracc de stefan}). Note that this kind of relationship, beside the strong maximum principle, has been a useful tool in the classical Stefan problems to analyse the existence and uniqueness of solutions. Finally we present an exact solution to a particular case of problem (\ref{St}) in terms of the Wright functions and the relationship already obtained is checked.

\section{Preliminaries}
\subsection{Basics of Fractional Calculus}

\begin{defi}\label{defi frac} Let $\left[a,b\right]\subset \bbR$ and $\al \in \bbR^+$ be such that $n-1<\al\leq n$. 
\begin{enumerate}
	\item   If $f \in L^1[a,b]$ we define the \textsl{fractional 
Riemann--Liouville integral of order  $\alpha$} as
$$_{a}I^{\alpha}f(t)=\frac{1}{\Gamma(\alpha)}\int^{t}_{a}(t-\tau)^{\alpha-1} f(\tau)d\tau. $$
\item If $f\in AC^n[a,b]=\left\{ f\, \in \,\CC^{(n-1)} \, | \, f^{(n-1)} \, \text{ is absolutly continuous} \right\}$, we define the 
\textsl{fractional Riemann--Liouville  derivative of order $\alpha$ } to
$$ ^{RL}_{a}D^{\alpha}f(t)= \left[ D^n_{a}I^{n-\alpha}f  \right] (t) =\frac{1}{\Gamma(n - \alpha)}\frac{d^n}{dt^n}\int^{t}_{a}(t-\tau)^{n-\alpha-1} f(\tau)d\tau. $$
\item If $ f \in W^n(a,b)=\left\{ f \, \in \, \CC^n(a,b]\, | \, f^{(n)}\in L^1[a,b] \right\}$ we define the \textsl{fractional Caputo derivative of order  $\al$} to
$$\,^C_{a} D^{\alpha}f(t)=\left\{\begin{array}{lc} \frac{1}{\Gamma(n-\al)}\displaystyle\int^{t}_{a}(t-\tau)^{n-\al-1} f^{(n)}(\tau)d\tau, &  n-1<\al<n\\
f^{(n)}(t), &   \al=n. \end{array}\right.$$
\end{enumerate}
\end{defi}

\begin{propo}\label{propo frac} The following properties involving the fractional integrals and derivatives hold:
\begin{enumerate}
\item \label{RL inv a izq de I} The  \textsl{fractional Riemann--Liouville derivative } operator is a left inverse of the \textsl{fractional Riemann--Liouville integral} of the same order  $\al\in \bbR^+$.\\ If $f \in L^1[a,b]$, then
$$^{RL}_{a}D^{\al}\,_{a}I^{\al}f(t)=f(t)  \quad a.e.$$

\item The fractional  Riemann--Liouville integral, in general,  is not a left inverse operator of the fractional derivative of Riemann--Liouville.\\
 
In particular, if $0<\al<1$, then 
$ _{a}I^{\al}(^{RL}_{a}D^{\al}f)(t)=f(t) - \dfrac{_{a}I^{1-\al}f(a^+)}{\G(\al)(t-a)^{1-\al}}.$ 

\item If $f\in AC^n[a,b],$ then   $^{RL}_{a}D^{\al}f (t)=\displaystyle\sum_{k=0}^{n-1} \frac{f^{(k)}(a)}{\G(1+k-\al)}(t-a)^{k-\al}+\, ^C_{a} D^{\alpha}f(t)$.
\end{enumerate}
\end{propo}

\begin{obs} It is known that if $f$ and $g$ are functions supported in $[0,\infty)$, then the convolution of $f$ and $g$ is defined by
$$(f\ast g)(t)=\int_0^tf(t-\tau)g(\tau)d\tau. $$
Then if we consider a function $f$ supported in $[0,\infty)$ and $\chi_\alpha$ is the locally integrable function defined by 
\begin{equation}\label{chi_alpha}\chi_\alpha(t)=\begin{cases} \frac{t^{\alpha-1}}{\G(\al)} & \text{if } t>0 \\
0 & \text{if } t\leq 0,\end{cases}
\end{equation}
 then we have the following properties:\\
\begin{equation}\label{I_al con convolucion}
 _{0}I^{\al}f(t)=\left(\chi_\alpha \ast f\right) (t)
\end{equation}
\begin{equation}\label{D_al RL con convolucion}
 _{0}^{RL}D^{\al}f(t)=\frac{d^n}{dt^n}\left(\chi_{n-\alpha} \ast f\right) (t)
\end{equation}
\begin{equation}\label{D_al C con convolucion}
 _{0}^CD^{\al}f(t)=(\chi_{n-\alpha}) \ast \frac{d^n}{dt^n}f (t)
\end{equation}

\end{obs}

\subsection{The special functions involved }

In this subsection we present some special functions that will be part of the explicit solution that will be presented in the next section. 

\begin{defi} For every $z\in \bbC$ , $\rho>-1$ and $\be\in \bbR$ the \textit{Wright} function is defined by
\begin{equation}\label{W} W(z;\rho;\be)=\sum^{\infty}_{k=0}\frac{z^{k}}{k!\G(\rho k+\be)}.\end{equation}
The  \textit{Mainardi} function $\cite{FM:1999}$ is a special case of the Wright function defined by  

\begin{equation}\label{M} M_\rho (z)= W(-z,-\rho,1-\rho)=\sum^{\infty}_{n=0}\frac{(-z)^n}{n! \G\left( -\rho n+ 1-\rho \right)}, \quad z\in \bbC, \, \rho<1.  \end{equation}

\end{defi}

\begin{obs} Function $(\ref{W})$ was presented by E. M. Wright at the beginning of the \textit{XX} Century, and he studied its asymptotic behaviour in $\cite{Wr1:1934}$ and  $\cite{Wr2:1940}$. It is known that: 
\begin{enumerate}
	\item  The Wright function is an entire function if $\rho>-1$.
	\item The derivative of the Wright function can be computed as
 \begin{equation}\label{derivada de W} \frac{\p}{\p z} W(z,\rho,\be) = W(z,\rho,\rho+\be).  \end{equation}

		\item Some particular cases are:\\
A Gaussian function: $\frac{1}{\sqrt{\pi}}e^{-x^2}= W\left(-2x,-\frac{1}{2},\frac{1}{2}\right)=M_{1/2}(2x).$\\
The error function: $ \mbox{erf\,}\left(x\right)=\frac{2}{\sqrt{\pi}}\displaystyle\int_{0}^{x} e^{-\xi^2}d\xi= 1-W\left(-2x,-\frac{1}{2},1\right)$. \\
And the complementary erf function: 
$$\mbox{erfc\,}\left(x\right)=\frac{2}{\sqrt{\pi}}\displaystyle\int_{x}^{\infty} e^{-\xi^2}d\xi=  
W\left(-2x,-\frac{1}{2},1\right).$$

\item \begin{equation}\label{lim_W=0} \displaystyle\lim_{x \rightarrow \infty } W\left(-x,-\frac{\al}{2},\beta\right)=0  \qquad \text{ for all } \, \al\in (0,1), \quad \beta>0.\end{equation}
\end{enumerate}
\end{obs}

\noindent The next two propositions were proved in \cite{RoSa1}.

\begin{propo}\label{M y W pos y decrec} If  $\,  0<\al<1$, then:
\begin{enumerate}
\item $M_{\al/2}(x) $ is a positive and strictly decreasing positive function in $\bbR^+$  such that $M_{\al/2}(x)<\frac{1}{\G\left(1-\frac{\al}{2}\right)}$ for all $x>0$.
\item  $W\left(-x,-\frac{\al}{2},1\right)$ is a positive and strictly decreasing function in  $\bbR^+$ such that  $0<W\left(-x,-\frac{\al}{2},1\right)< 1, $ for all $x>0$.
\item  $1-W\left(-x,-\frac{\al}{2},1\right)$ is a positive and strictly increasing function in  $\bbR^+$ such that  $0<1-W\left(-x,-\frac{\al}{2},1\right)< 1, $ for all $x>0$.
\end{enumerate}
\end{propo}

\begin{propo}\label{conv M y W cuando al tiende a 1}If $x\in \bbR^+_0$ and $\al \in (0,1)$ then:
\begin{enumerate}
\item $\lim\limits_{\al\nearrow 1}M_{\al/2}\left(2x\right)=M_{1/2}(2x)=\frac{e^{-x^2}}{\sqrt{\pi}}$;
\item $\lim\limits_{\al\nearrow 1}\left[1-W\left(-2x,-\frac{\al}{2},1\right)\right]= erf\left(x\right) .$
\end{enumerate}
\end{propo}

\section{The fractional Stefan Problem}

Let us now study problem (\ref{St}), which is the main goal of this paper.  The following two regions will be considered: $\Omega_T=\{ (x,t) / 0<x<s (t), \, 0<t\leq T \}$ and   
 $\p_p \Omega_T=\{ (0,t), 0<t\leq T \} \cup \{ (s(t),t), 0<t \leq T \} \cup \{(x,0), 0\leq x \leq b\}$, where the latter is called parabolic boundary.

\begin{defi}\label{Def sol St} A pair $\{u,s\} $ is a solution of  problem $(\ref{St})$  if

\begin{enumerate}
    \item $u$ is defined in $[0,b_0]\times[0,T]$ where $b_0:= \max\{s(t), 0 \leq t \leq T \}$.
	\item  $u\in $ $C(D_T)\cap C^2_x(D_T)$, such that $u_x \in 
	AC^1_t((0,T))$ 	where \\$AC^1_t((0,T)):=\{f(x,\cdot)\colon f \in AC^1(0,T)  \quad \text{for every fixed } x\in [0,b_0] \}$.
\item $u$ is continuous in $D_T \cup \p_p D_T$ except  perhaps at $(0,0)$ and $(b,0)$ where  
 $$ 0\leq \underset{(x,t)\rightarrow (0,0)}{\liminf}u(x,t)\leq \underset{(x,t)\rightarrow (0,0)}{\limsup } u(x,t)<+\infty$$
  and 
	$$0\leq \underset{(x,t)\rightarrow (b,0)}{\liminf}u(x,t)\leq \underset{(x,t)\rightarrow (b,0)}{\limsup } u(x,t)<+\infty. $$ 
    \item $s \in C^1(0,T)$.
		\item There exists $\left.^{RL}_{0}D^{1-\al}_t \frac{\p}{\p x} u (x,t)\right|_{(s(t),t)}$ for all $t \in (0,T]$.
    \item $u$ and $s$ satisfy $(\ref{St})$.
   	
		\end{enumerate}
\end{defi}

\begin{obs}\label{C->RL}  Similar problems for the fractional derivative in the Caputo sense were studied in  \cite{Atkinson,  Ro:2016,  RoSa1, RoSa2, RoTa:2014, VMD:2012, Voller:2014, Voller:2013}. For example, the formulation given in \cite{RoSa1}  is: 
\begin{equation}{\label{St-Caputo}}
\begin{array}{llll}
     (i)  &  ^CD_t^{\al} u(x,t)=\lambda^2 \frac{\p}{\p x^2}u(x,t), &   0<x<s(t), \,  0<t<T, \,  0<\al<1 , \, \, \la>0\\
     (ii) &  u(x,0)=f(x), & 0\leq x\leq b=s(0),\\ 
       (iii) &  u(0,t)=g(t),   &  0<t\leq T,  \\
         (iv) & u(s(t),t)=0, & 0<t\leq T,\\
  (v) & ^CD^{\al}s(t)=-u_x(s(t),t), & 0<t\leq T.                                              \end{array}\end{equation}

We assert that problem (\ref{St-Caputo}) is similar to problem (\ref{St}) because from Proposition 1 it results that   if $u$ is a solution  of the fractional diffusion equation  for the Caputo derivative  $(\ref{St-Caputo}-i)$, then $u$ verifies the FDE $(\ref{St}-i)$.\\

However the converse of the previous statement is not true because the fractional integral of Riemann-Liouville is not the inverse operator of the Riemann-Liouville derivative of equal order. 

Also is worth noting that if we apply the integral operator $^{RL}_{0}D^{1-\al}_t$ to both members of equation $(\ref{St-Caputo}-v)$, we get equation 
$$ \frac{d}{dt}s(t)=-^{RL}_{0}D^{1-\al}_t \left[\frac{\p}{\p x} u (s(t),t)\right], \quad 0<t\leq T,  $$
which is different to equation (\ref{cond fracc de stefan}) unless $\al=1$.  
\end{obs}

\subsection{An exact solution}
Without loss of generality, we take now $\lambda=1$. Also we consider  $b=0$ and  a constant boundary condition in $(\ref{St}-iii)$. Namely, let the following fractional Stefan problem be:
\begin{equation}{\label{St_1}}
\begin{array}{llll}
     (i)  &   \frac{\p}{\p t}u(x,t)=\, \frac{\p}{\p x}\left(^{RL}_{0}D^{1-\al}_t \frac{\p}{\p x} u (x,t)\right), &   0<x<s(t), \,  0<t<T, \,  0<\al<1 , \, \,\\
 
       (ii)  &  u(0,t)=1,   &  0<t\leq T,  \\

         (iii) & u(s(t),t)=0, & 0<t\leq T, \, s(0)=0 \\
        
       (iv) & \frac{d}{dt}s(t)=-\left.^{RL}_{0}D^{1-\al}_t \frac{\p}{\p x} u (x,t)\right|_{(s(t),t)}, & 0<t\leq T, 
                                             \end{array}\end{equation}

\noindent In order to find an exact solution to problem (\ref{St_1}), the next lemma proved in \cite{GoReRoSa:2015} will be used:
\begin{lema}\label{int_x c es sol} Let $c(x,t)$ be a solution of the time--fractional diffusion equation for the Caputo derivative  $(\ref{St-Caputo}-i)$ such that:
\begin{equation}\label{hip lema integral 1}
\quad \text{ For every } (x,t),\text{ the function } F(x,t)=\int^{\infty}_{x}c(\xi,t)d\xi  \text{ is well defined, }
  \end{equation}
 \begin{equation}\label{hip lema integral 2}
\quad \displaystyle\lim_{x \rightarrow  \infty}\dfrac{\partial c}{\partial x}(x,t)=0,\,
 \end{equation}
\begin{equation}\label{hip lema integral 3}
 \quad  \left|\frac{\partial}{\partial \tau}c(\xi,\tau)\right|\leq g(\xi) \in L^1(x,\infty)\, , 
\end{equation}
\begin{equation}\label{hip lema integral 4}
  \quad \frac{\frac{\partial}{\partial \tau}c(\xi,\tau)}{(t-\tau)^\al} \in L^1((x,\infty)\times(0,t))\, . 
 \end{equation}

\noindent Then $\displaystyle\int^{\infty}_{x}c(\xi,t)d\xi$ is a solution to the time fractional diffusion equation for the Caputo derivative  $(\ref{St-Caputo}-i)$.

\end{lema}

\begin{obs} The factor $2$ appearing in the next functions $(\ref{u})$ and $(\ref{s})$ was considered with the aim to recover the Gaussian and erf functions when we make $\al \nearrow 1$ (according to Proposition $\ref{conv M y W cuando al tiende a 1}$).

\end{obs}

\begin{teo}\label{teo exact sol} The pair given by 
\begin{equation}{\label{u}} u(x,t)=1-\frac{1}{1-W\left(-2\xi,-\frac{\al}{2},1\right)}[1-W\left(-\frac{x}{t^{\al/2}},-\frac{\al}{2},1\right)],
\end{equation}
\begin{equation}{\label{s}} s(t)=2\xi t^{\al/2},
\end{equation}
where $\xi$ is the unique solution to the equation 
\begin{equation}{\label{eq_xi}} 2x\left[ 1-W\left(-2x,-\frac{\al}{2},1\right)\right]=2xW\left(-2x,-\frac{\al}{2},1\right)+ W\left(-2x,-\frac{\al}{2},1+\frac{\al}{2}\right), \quad x>0
\end{equation}
is a solution to problem \rm{(\ref{St_1})}.

\end{teo}
\proof We know that 
\begin{equation} u(x,t)=a+b \left[1-W\left(-\frac{x}{ t^{\al/2}},-\frac{\al}{2},1\right)\right] \end{equation}
 is a  solution of fractional diffusion equation for the Caputo derivative ($\ref{St-Caputo}-i)$ for all $a \in \bbR$, $b \in \bbR$ (see \cite{FM-TheFundamentalSolution} or \cite{RoTa:2014}). 
Then, applying Remark \ref{C->RL}, it results that $u$ is a solution to equation $(\ref{St_1}-i)$.\\
From $(\ref{St_1}-ii)$ we obtain
\begin{equation}\label{St-ii} 1=u(0,t)=a+b \left[1-W\left(0,-\frac{\al}{2},1\right)\right]=a,
\end{equation}
and from $(\ref{St_1}-iii)$ we get 
\begin{equation}\label{St_1-iii} u(s(t),t)=1+b\left[1- W\left(-\frac{s(t)}{ t^{\al/2}},-\frac{\al}{2},1\right)\right]=0.
\end{equation}
 Note that (\ref{St_1-iii}) must be verified for all $t>0$, so we will ask for  $s(t)$ to be proportional to $t^{\al/2}$, that is to say
\begin{equation}\label{s(t)}
s(t)= 2\xi t^{\al/2} \qquad  \text{for some } \  \xi>0.
\end{equation}
 Replacing (\ref{s(t)}) in (\ref{St_1-iii})  and taking into account Proposition \ref{M y W pos y decrec}  it follows that   
$b= -\frac{1}{1-W\left(-2\xi,-\frac{\al}{2},1\right)}$ and then $(\ref{u})$ holds.
 
With the aim of use the fractional Stefan condition $(\ref{St_1}-iv)$, the fractional derivative  $ ^{RL}_{0}D^{1-\al}_t \left(\frac{\p}{\p x} u (x,t)\right)$ must be computed.\\

From  (\ref{lim_W=0}) and using estimates made in \cite{GoReRoSa:2015}, it yields that, for every $x>0$, $w_1(x,t)=W\left(-\frac{x}{t^{\al/2}}, -\frac{\al}{2},1\right)$ is under the assumptions of Lemma \ref{int_x c es sol}.  Clarely, $w_2(x,t)=x$ is a solution to the FDE $(\ref{St_1}-i)$.

Then, using the linearity of the Caputo derivative \cite{Kilbas} and the principle of superposition we can state that the function defined by
\begin{equation}{\label{2} }
v(x,t)= -\left[1-\frac{1}{1-W\left(-2\xi,-\frac{\al}{2},1\right)}\right]x+\frac{t^{\al/2}}{1-W\left(-2\xi,-\frac{\al}{2},1\right)}W\left(-\frac{x}{t^{\al/2}},-\frac{\al}{2},1+\frac{\al}{2}\right) \end{equation}
is a solution of the FDE such that  $\frac{\p v}{\p x}(x,t)=-u(x,t)$ for all $x>0, t>0$. Hence
$$ \hspace{-6.5cm}\frac{\p}{\p t} v(x,t)= \, \frac{\p}{\p x}\left(^{RL}_{0}D^{1-\al}_t \frac{\p}{\p x} v (x,t)\right)=$$
\begin{equation}{\label{1} }
=\, \frac{\p}{\p x}\left(^{RL}_{0}D^{1-\al}_t (-u(x,t))\right)=-\,^{RL}_{0}D^{1-\al}_t \frac{\p}{\p x} u (x,t), \, x>0, \, t>0. \end{equation}
Derivating $v$  with respect to the $t$ variable and using (\ref{1}) it yields 
\begin{equation}\label{D1-al_u_x}
-\,^{RL}_{0}D^{1-\al}_t \frac{\p}{\p x} u (x,t)=\frac{1}{1-W\left(-2\xi,-\frac{\al}{2},1\right)}\frac{\al}{2}
\left[\frac{x}{t}W\left(-\frac{x}{ t^{\al/2}},-\frac{\al}{2},1\right)
+t^{\al/2-1}W\left(-\frac{x}{ t^{\al/2}},-\frac{\al}{2},1+\frac{\al}{2}\right)\right].
\end{equation}

Replacing $(\ref{s(t)})$ and $(\ref{D1-al_u_x})$ into the fractional Stefan condition $(\ref{St_1}-iv)$ it results that $\xi$  must verify the following equation:
\begin{equation}\label{eq xi}
2 x \left[ 1-W\left(-2x,-\frac{\al}{2},1\right) \right]= 2 x W\left(-2x,-\frac{\al}{2},1\right)
+ W\left(-2x,-\frac{\al}{2},1+\frac{\al}{2}\right).
 \end{equation}

Define functions $H(x)=x \left[ 1-W\left(- x,-\frac{\al}{2},1\right) \right]$  and $G(x)= x W\left(-x,-\frac{\al}{2},1\right)
+ W\left(-x,-\frac{\al}{2},1+\frac{\al}{2}\right)$ in $\bbR_0$.
From Proposition \ref{M y W pos y decrec}, H is an increasing function such that $H(0)=0$. On the other hand, $G$ is a decreasing function in $\bbR^+$ such that 
$G(0)=\frac{1}{\G \left(1+\frac{\al}{2}\right)} > 0 $ due to  $G '(x)=-xM_{\al/2}(x)<0$ for all $x>0$ and $\al \in (0,1)$. Then, we can assert that there exists a unique positive solution $\xi$ such that $H(2\xi)=G(2\xi)$.

\endproof

\subsection{An integral relationship between $u$ and $s$}
 
The next Theorem provides an integral relationship between the free boundary $s$ and function $u$, obtained from  the fractional Stefan condition (\ref{cond fracc de stefan}).

\begin{teo}\label{Teo cond de Stefan integral}
Let  $\{u,s\}$ be a solution of problem  $(\ref{St})$ such that $\frac{\p^2}{\p t\p x} u (x,t) \in \CC^1(\Omega_T) $, $g \in \, AC^{1}(0,T)$ and,  $^{RL}D_t^{1-\al}g$ and $^{RL}D_t^{1-\al} u(x,t)|_{(s(t),t)} \in L^1(a,b)$. Then  the following integral condition for the free boundary $s(t)$ and the function $u(x,t)$
\begin{equation}\label{cond de Stefan integral}
 s^2(t)=b^2+2\int_0^t\,^{RL}D_t^{1-\al}g(\tau)d\tau + 2\int_0^b z f(z)dz -2\int_0^{s(t)}z u(z,t)dz-2\int_0^t\,\left. ^{RL}D_t^{1-\al}u(x,t)\right|_{(s(\tau),\tau)}d\tau 
\end{equation}
is verified.
\end{teo}

\proof
Recall the Green identity: 
$$ \int_{\p \Omega}Pdt+Qdx =\iint_{\Omega} (Q_t-P_x )\, dA , $$
where $\Omega$ is an open simply connected region,  $\p \Omega$ is a positively oriented, piecewise smooth, simple closed curve, and the field $F=(P,Q)$ is $\CC^1$ in an open set containing $\Omega$.\\

\noindent Consider the functions $P$ and $Q$ defined by
 \begin{equation}\label{P} P(x,t)=-x\,^{RL}D_t^{1-\al} u_x(x,t) +\, ^{RL}D_t^{1-\al} u(x,t)\end{equation}  
\begin{equation}\label{Q} Q(x,t)=-x\, u(x,t).\end{equation}  
and the region  $$\Omega_\e=\left\{(x,\tau) \in \bbR^2 \, / \, \e<\tau<t, 0<x<s(\tau)\right\}, \, \e >0.$$
\noindent Note that 
$$^{RL}D_t^{1-\al} u_x(x,t)= \frac{\p }{\p t} I^{1-\al}_t(u_x(x,t))= \frac{\p }{\p t}\left( \chi_{1-\al} (t)\ast u_x(x,t)\right),$$
where $\chi_\al$ (defined in eq. (\ref{chi_alpha}) ) is an $L^1_{loc}$ function. Then the convolution  inherits all the properties of $u_x$ and, taking into account that $u$ is under the assumptions of Definition \ref{Def sol St} and that $\frac{\p^2}{\p t\p x}u(x,t) \in \CC^1(\Omega_T) $, the field $F$ have all the regularity required in  $\Omega_\e$. \\

\noindent Also, due to the regularity of the field F, the derivatives $\p/ \p x$  and $\,^{RL}D_t^{1-\al}$ commutes. Then, applying Green's Theorem  and taking into account that $u$ verifies $(\ref{St}-i)$, we get 
  
$$\hspace{-5cm} \int_{\partial \Omega_\e}Pd\tau+Qdx=\int_{\p \Omega_\e} \left[-x\,^{RL}D_t^{1-\al} u_x(x,\tau) +\, ^{RL}D_t^{1-\al} u(x,\tau)\right] d\tau-x\, u(x,t) dx = $$
$$\hspace{2cm} =\iint_{\Omega_\e}\left[ -x u_t(x,\tau)+
\,^{RL}D_t^{1-\al} u_x(x,\tau)+x \frac{\p}{\p x}\left( \,^{RL}D_t^{1-\al} u_x(x,\tau)  \right)-\frac{\p}{\p x} \left(\,^{RL}D_t^{1-\al} u(x,\tau)\right) \right] d\tau\, dx $$
\begin{equation}\label{int de linea cero}\hspace{-5cm}=\iint_{\Omega_\e} x \left[\,^{RL}D_t^{1-\al} u_{xx}(x,\tau)- u_t(x,\tau)\right] d\tau\, dx = 0.\end{equation}
%

\noindent Consider
$ \p \Omega_\e= \p \Omega_{\e1}\cup \p \Omega_{\e2}\cup \p \Omega_{\e3}\cup \p \Omega_{\e4} $ where
 $\p \Omega_{\e1}= \left\{(0,\tau) , \, \e \leq \tau \leq t \right\}$,
$\p \Omega_{\e2}= \left\{(z,\e) , \, 0 \leq z \leq s(\e) \right\}$,
 $\p \Omega_{\e3}= \left\{(s(\tau),\tau) , \, \e \leq \tau \leq t \right\}$ and  
$\p \Omega_{\e4}= \left\{(z,t) , \, 0 \leq z \leq s(t) \right\}$.

Integrating over the contour  $\p \Omega_\e$ (positively oriented) we get: 
\begin{equation}\label{int_pD_e1}
\int_{\p \Omega_{\e1}} P d\tau+ Q dx =\int^\e_t\, ^{RL}D_t^{1-\al} u(0,\tau)d\tau=-\int_\e^t\, ^{RL}D_t^{1-\al} g(\tau) d\tau,
\end{equation}
\begin{equation}\label{int_pD_e2} 
\int_{\p \Omega_{\e2}} P d\tau+ Q dx =\int_0^{s(\e)}-z u(z,\e)dz,
\end{equation}
\begin{equation}\label{int_pD_e3-antes-de-reempl}
\hspace{-0.5cm}\begin{array}{rcl}
 \displaystyle\int_{\p \Omega_{\e3}} P d\tau+ Q dx  & = &   \displaystyle\int_\e^t \left[-s(\tau)\,\left.^{RL}_{0}D^{1-\al}_t \frac{\p}{\p x} u (x,t)\right|_{(s(\tau),\tau)} +\right.
 \end{array}
\end{equation}
$$\hspace{4cm} \left. \,^{RL}D_t^{1-\al} \left. u (x,t)\right|_{(s(\tau),\tau)} - s(\tau)u(s(\tau),\tau) s'(\tau)\right] d\tau $$
\begin{equation}\label{int_pD_e3}
\hspace{2.8cm}\begin{array}{rcl}
& = & \displaystyle\int_\e^t s(\tau)s'(\tau) d\tau+ \displaystyle\int_\e^t\,^{RL}D_t^{1-\al} \left. u (x,t)\right|_{(s(\tau),\tau)} d\tau\\
 & = & \dfrac{s^2(t)}{2}-\dfrac{s^2(\e)}{2} + \displaystyle\int_\e^t\,^{RL}D_t^{1-\al} \left. u (x,t)\right|_{(s(\tau),\tau)} d\tau,
\end{array}
\end{equation}
\begin{equation}\label{int_pD_e4}
\int_{\p \Omega_{\e4}} P dt+ Q dx =\int_0^{s(t)}z u(z,t)dz.
\end{equation}

\noindent Join (\ref{int de linea cero}), (\ref{int_pD_e1}), (\ref{int_pD_e2}), (\ref{int_pD_e3}) and (\ref{int_pD_e4}), it results that 

\begin{equation}
-\int_\e^t\, ^{RL}D_t^{1-\al} g(\tau) d\tau-\int_0^{s(\e)}z u(z,\e)dz +\frac{s^2(t)}{2}-\frac{s^2(\e)}{2}+ \displaystyle\int_\e^t\,^{RL}D_t^{1-\al} \left. u (x,t)\right|_{(s(\tau),\tau)} d\tau+ 
\end{equation}
$$ \hspace{9cm}+\int_0^{s(t)}z u(z,t)dz=0. $$
\noindent Taking the limit when $\e \searrow 0$ we get the integral relationship (\ref{cond de Stefan integral}), i.e. the thesis holds.

\endproof

\begin{obs} If we take $\al=1$ in the integral relationship (\ref{cond de Stefan integral}) we get
$$ s^2(t)=b^2+2\int_0^t\,g(\tau)d\tau + 2\int_0^b z f(z)dz -2\int_0^{s(t)}z u(z,t)dz-2\int_0^t\, u(s(\tau),\tau)d\tau,$$ 
and using condition $(\ref{St}-iii)$, it results that 
\begin{equation}\label{cond de Stefan integral clasica}
s^2(t)=b^2+2\int_0^t\,g(\tau)d\tau + 2\int_0^b z f(z)dz -2\int_0^{s(t)}z u(z,t)dz, 
\end{equation}
 where (\ref{cond de Stefan integral clasica}) is the classical integral condition for the free boundary when a classical Stefan problem is considered (see\cite{Cannon}--Lemma 17.1.1). \\
It was also   proven in \cite{Cannon} that condition  (\ref{cond de Stefan integral clasica}) is equivalent to the Stefan condition
\begin{equation}\label{cond de Stefan clasica}
\frac{d}{dt}s(t)=-\frac{\p }{\p x}u(s(t), t), \quad \forall \, t>0. 
\end{equation} 
 Hence, it is natural to wonder if the ``fractional Stefan condition'' (\ref{cond fracc de stefan}) and the ``fractional integral condition'' (\ref{cond de Stefan integral}) are equivalent too.
\end{obs} 

\begin{teo}\label{reciproco - cond de Stefan integral}
Let  $\{u,s\}$ be a solution of problem   $\left\{ (\ref{St}-i), (\ref{St}-ii), (\ref{St}-iii), (\ref{St}-iv), (\ref{cond de Stefan integral})\right\}$ such that $\frac{\p^2}{\p t\p x} u (x,t) \in \CC^1(\Omega_T) $,   $g \in \, AC^{1}(0,T)$,  $^{RL}D_t^{1-\al}g$,   $^{RL}D_t^{1-\al} u(x,t)|_{(s(t),t)} \in L^1(a,b)$ and $s(t)>0$ for all $t>0$. Then functions  $s=s(t)$ and  $u=u(x,t)$ verify the fractional Stefan condition (\ref{cond fracc de stefan}) .
\end{teo}
\proof

Reasoning as in Theorem \ref{Teo cond de Stefan integral}, we can state that the equalities (\ref{int de linea cero}), (\ref{int_pD_e1}),  (\ref{int_pD_e2}) and (\ref{int_pD_e3-antes-de-reempl}) hold. Then, taking the limit when $\e \searrow 0$ it results that
$$
-\int_0^t\, \hspace{-0.2cm}^{RL}D_t^{1-\al} g(\tau) d\tau-\int_0^{b}\hspace{-0.4cm}z f(z)dz -\displaystyle\int_0^t\,s(\tau) ^{RL}D_t^{1-\al} \left. \frac{\p}{\p x}u (x,t)\right|_{(s(\tau),\tau)} d\tau + $$
\begin{equation}\label{teo 2-1}
 \hspace{3cm} +\displaystyle\int_0^t\, ^{RL}D_t^{1-\al} \left. u (x,t)\right|_{(s(\tau),\tau)} d\tau + \int_0^{s(t)}\hspace{-0.4cm}z u(z,t)dz=0. 
\end{equation}

Multiplying (\ref{teo 2-1}) by 2 and using hypothesis (\ref{cond de Stefan integral}) it yields that 
\begin{equation}\label{teo 2-2}
-2\displaystyle\int_0^t\,s(\tau) ^{RL}D_t^{1-\al} \left. \frac{\p}{\p x}u (x,t)\right|_{(s(\tau),\tau)} d\tau = s^2(t)-b^2.
\end{equation}
 Differentiating both sides of equation (\ref{teo 2-2}) whith respect to the $t-$variable and being $s(t)>0$ for all $t>0$, the thesis holds.

\endproof

\begin{obs} The hypothesis $s(t)>0$ for all $t>0$ in the previous Theorem is not necessary in the classical Stefan problem. In fact, the Stefan condition  (\ref{cond de Stefan clasica})  join with the maximum principle imply that $u$ is a decreasing function of the $x-$variable for every $t>0$,  leading function $s$ to be a non-decreasing function of $t$.
\\
However this simple tool can not be considered in this case, because decreasing functions may have a positive Riemann-Liouvulle derivative. For example, let   $\al \in (0,1)$ be and consider $\gamma \in (0,1) $ such that $0<\al - \gamma$. Function $f(t)=t^{-\gamma}$, $t>0$ is a decreasing function in $\bbR^+$ and 
$$ ^{RL}D^{1-\al} (t^{-\gamma})=\frac{\G(-\gamma +1)}{\G(-\gamma -(1-\al) +1)}t^{-\gamma+\al-1}= \frac{\G(-\gamma +1)}{\G(\al-\gamma)}t^{-\gamma+\al-1}>0 \,$$
for all $t>0$.

\end{obs}

\subsection{Example}
 The solution (\ref{u}), (\ref{s}) to problem (\ref{St_1}) given in Theorem \ref{teo exact sol} verify the integral relationship (\ref{cond de Stefan integral}). In this case, $g(t)=1$ for all $t>0$, $b=0$, $s(t)=2\xi t^{\al/2}$ and  $u(x,t)=1-\frac{1}{1-W\left(-2\xi,-\frac{\al}{2},1\right)}[1-W\left(-\frac{x}{t^{\al/2}},-\frac{\al}{2},1\right)]$
where $\xi$ is the unique solution to equation (\ref{eq_xi}).\\
The fractional derivative of Riemann-Liouville of a constant is easy to compute (see \cite{Podlubny}) and it is given by
\begin{equation}\label{check-1}
^{RL}D_t^{1-\al}1=\frac{t^{0-(1-\al)}}{\G (0-(1-\al)+1 )}=\frac{t^{\al-1}}{\G (\al) }.
\end{equation}
Integrating $(\ref{check-1})$ from $0$ to $t$, and using the Gamma function property $z\G(z)=\G(z+1)$, we have 
\begin{equation}\label{check-2}
\int_0^t\,^{RL}D_t^{1-\al}1d\tau=\frac{t^{\al}}{\G(\al+1)}.
\end{equation}

\noindent Applying  two times Lemma \ref{int_x c es sol} to function $W\left(-\frac{x}{t^{\al/2}},-\frac{\al}{2},1\right)$, we get that function 
$w(x,t)=W\left(-\frac{x}{t^{\al/2}},-\frac{\al}{2},1+\al\right)t^\al$ is a solution to the fractional diffusion equation for the Caputo derivative $(\ref{St-Caputo}-i)$ and therefore is a solution to $(\ref{St}-i)$ such that $\frac{\p^2}{\p x^2}w(x,t)=
W\left(-\frac{x}{t^{\al/2}},-\frac{\al}{2},1\right)$. Then  
\begin{equation}\label{check-3}
\frac{\p}{\p t}w(x,t)=\,^{RL}D_t^{1-\al}\frac{\p^2}{\p x^2}w(x,t)=^{RL}D_t^{1-\al}W\left(-\frac{x}{t^{\al/2}},-\frac{\al}{2},1\right).
\end{equation}
Using (\ref{check-3}) and the linearity of the Riemann-Liouville derivative, it results that

$$\hspace{-10cm}\,^{RL}D_t^{1-\al} u(x,t)=$$
$$=\,^{RL}D_t^{1-\al}\left(1-\frac{1}{1-W\left(-2\xi,-\frac{\al}{2},1\right)}\right)+\frac{1}{1-W\left(-2\xi,-\frac{\al}{2},1\right)} \frac{\p}{\p t}w(x,t).$$
Hence 
$$\hspace{-2.5cm}
\left.\,^{RL}D_t^{1-\al} u(x,t)\right|_{(2\xi t^{\al/2},t)}=
\left(1-\frac{1}{1-W\left(-2\xi,-\frac{\al}{2},1\right)}\right)\frac{t^{\al-1}}{\G(\al)}+$$
\begin{equation}\label{check-4}+\frac{t^{\al-1}}{1-W\left(-2\xi,-\frac{\al}{2},1\right)} \left[
W\left(-2\xi,-\frac{\al}{2},1+\frac{\al}{2}\right)\al \xi +W\left(-2\xi,-\frac{\al}{2},1+\al\right)\al \right].
\end{equation}
Integrating (\ref{check-4}) from  $0$ to $t$, 

$$\hspace{-2cm} \int_0^t\left.\,^{RL}D_t^{1-\al} u(x,t)\right|_{(2\xi \tau^{\al/2},\tau)}d\tau = \left(1-\frac{1}{1-W\left(-2\xi,-\frac{\al}{2},1\right)}\right)\frac{t^\al}{\G (\al+1)}+$$
\begin{equation}\label{check-5}\hspace{1cm}+\frac{ t^{\al}}{1-W\left(-2\xi,-\frac{\al}{2},1\right)} \left[
W\left(-2\xi,-\frac{\al}{2},1+\frac{\al}{2}\right) \xi +W\left(-2\xi,-\frac{\al}{2},1+\al\right)\right].
\end{equation}

\noindent Integrating by parts, the next computation follows:

$$\hspace{-10cm}
\int_0^{2\xi t^{\al/2}}z u(z,t)dz= $$
$$=2\xi^2t^\al\left(1-\frac{1}{1-W\left(-2\xi,-\frac{\al}{2},1\right)} \right)-\frac{1}{1-W\left(-2\xi,-\frac{\al}{2},1\right)}W\left(-2\xi,-\frac{\al}{2},1+\al\right) t^{\al}+ 
$$
\begin{equation}\label{check-6}
-\frac{2\xi t^{\al}}{1-W\left(-2\xi,-\frac{\al}{2},1\right)}
W\left(-2\xi,-\frac{\al}{2},1+\frac{\al}{2}\right)+\frac{ t^\al}{\left(1-W\left(-2\xi,-\frac{\al}{2},1\right)\right)\G(1+\al)}.
\end{equation}

\noindent Taking  account (\ref{check-2}), (\ref{check-5}) and (\ref{check-6}), it follows that
$$
2\int_0^t\,^{RL}D_t^{1-\al}1d\tau-2\int_0^t\left.\,^{RL}D_t^{1-\al} u(x,t)\right|_{(2\xi \tau^{\al/2},\tau)}-2\int_0^{2\xi t^{\al/2}}z u(z,t)dz= $$
%
%
$$ =-4\xi^2t^\al + \frac{4\xi^2t^\al}{1-W\left(-2\xi,-\frac{\al}{2},1\right)} + \frac{2\xi t^\al}{1-W\left(-2\xi,-\frac{\al}{2},1\right)}W\left(-2\xi,-\frac{\al}{2},1+\frac{\al}{2}\right) $$

\begin{equation}\label{check-7}
=2\xi t^\al\left[-2\xi +\frac{2\xi}{1-W\left(-2\xi,-\frac{\al}{2},1\right)} +\frac{1 }{1-W\left(-2\xi,-\frac{\al}{2},1\right)}W\left(-2\xi,-\frac{\al}{2},1+\frac{\al}{2}\right)\right].
\end{equation}
Finally, taking into account that $\xi$ verify equation $(\ref{eq_xi})$, it results that
$$
2\int_0^t\,^{RL}D_t^{1-\al}1d\tau-2\int_0^t\left.\,^{RL}D_t^{1-\al} u(x,t)\right|_{(2\xi \tau^{\al/2},\tau)}-2\int_0^{2\xi t^{\al/2}}\hspace{-0.8cm}z u(z,t)dz=4\xi^2t^\al=s^2(t), $$
and the fractional integral condition $(\ref{cond de Stefan integral})$ is satisfied. 

\section{Conclusions}
We have defined a new fractional Stefan problem (\ref{St}) by using the fractional Riemann-Liouville derivative. We have obtained an integral relationship between the fractional temperature and the fractional free boundary, which is equivalent to the fractional Stefan condition. We have also shown an exact solution to the fractional Stefan problem   (\ref{St_1}), which verifies the integral relationship $(\ref{cond de Stefan integral})$.

\section*{Acknowledgments}
The present work has been sponsored by the Projects PIP N$^\circ$ 0534 from CONICET--Univ. Austral, and by AFOSR--SOARD Grant FA$9550-14-1-0122$.

\bibliographystyle{plain}

\bibliography{RoTaBIBLIO}

\end{document}